\documentclass[a4paper,11pt,draft]{amsart}

\renewcommand{\Bbb}{\mathbb}
\newcommand{\TT}{{\mathbb T}}
\newcommand{\cal}{\mathcal}
\theoremstyle{plain}
\newtheorem{Th}{Theorem}
\newtheorem{Prop}{Proposition}
\newtheorem{Cor}{Corollary}
\newtheorem{Lem}{Lemma}
\newtheorem{Conj}{Conjecture}

\numberwithin{equation}{section}

\begin{document}

%Topmatter
\title [Strongly exposed points in the ball of the Bergman space]{
Strongly exposed points in the ball of the Bergman space}
\author{Paul Beneker}
\address{University of Amsterdam\\
         Korteweg de Vries Institute of Mathematics\\
         Plantage Muidergracht 24\\
         1018 TV Amsterdam\\
         The Netherlands}
\email[Paul Beneker]{beneker@@science.uva.nl}
\thanks{Research of the first author supported by the
        Netherlands research organization NWO}
\author{Jan Wiegerinck}
\email[Jan Wiegerinck]{janwieg@@science.uva.nl}
\date{July 3, 2002}

%End Topmatter
\begin{abstract}
    We investigate which boundary points in the closed unit ball of the
    Bergman space $A^1$ are strongly exposed. This requires study of the 
    Bergman projection and its kernel, the annihilator of Bergman space. 
    We show that all polynomials in the boundary of the unit ball are
    strongly exposed. 

\end{abstract}
\maketitle

\section{Introduction.} \label{S:intro}
In Banach space theory one often seeks to determine the geometry
of the unit ball of a given Banach space. A common way to distinguish
``round" and ``flat" parts of the boundary of the unit ball is through
extreme and non-extreme points. Among the {\sl extreme} points, or ``round"
parts of the boundary, further refinements can be made, for example
{\sl exposed} and {\sl strongly exposed} points. In this paper we 
study these sets for the (unweighted) Bergman space $A^1$ of the unit disc $D$
of $\Bbb C$.
These questions were inspired by and can be stated in the context of
Hardy spaces of the unit ball in ${\Bbb C}^n$. However, we attempt to
frame them within the theory of Bergman spaces. For an excellent survey 
of the theory of Bergman spaces, we refer the reader to ~\cite{hkz}. Our 
main result identifies a large class of strongly exposed points, which
includes all normalized polynomials. We also exhibit exposed points which
are not strongly exposed. In the process we find opportunity to study the 
subspace $(A^1)^\perp+C(\overline D)$ of $L^\infty(D)$, which is the analog of
$H^\infty+C\subset L^\infty(\TT)$ ($\TT=\partial D$).
As is $H^\infty+C$, the space $(A^1)^\perp+C(\overline D)$ is closed,
but contrary to
 $H^\infty+C$ it turns out not to be an algebra. However, it is a $C$-module.

The authors thank M.A.~Kaashoek  and R.~Kortram for stimulating discussions.

\section{The Bergman Space $A^1$.} \label{S:Bergman}
Let $D$ be the (open) unit disc in $\Bbb C$, with boundary
$\Bbb T$, and let $dA=\frac{1}{\pi}r\,dr\,d\theta$ be the normalized
Lebesgue measure on $D$. The space of all holomorphic functions on
$D$ will be denoted by $H(D)$. The (unweighted) {\sl Bergman space} $A^1$
consists of all functions $f$ in $H(D)$ that are area-integrable on
$D$. Equipped with the {\sl Bergman norm},
$$
  \| f\| := \int_D |f(z)| \,dA(z),
$$
the Bergman space becomes a Banach space.

More generally, for $1\le p < \infty$, the space $A^p =H(D) \cap
L^p(D,dA)$ with the $L^p$-norm is a Banach space. Under the $L^2$-inner
product, $A^2$ is a Hilbert space. The orthogonal projection
$P:L^2\to A^2$, the so-called Bergman projection, will play an important role 
later on, so we mention its explicit representation:
\begin{eqnarray}\label{eq:bergP}
   Pf (z) &=& = \int_D \frac{f(w)}{(1-z\overline{w})^2} \,dA(w)\label{eq:Bproj}\\
          &=& \sum_{n=0}^\infty \bigl((n+1)\cdot\int_D f(w)\overline{w}^n \,dA(w)\bigr) z^n\label{eq:bergP2}.
\end{eqnarray}

\section{Extreme, exposed and strongly exposed points.} \label{S:expos}
Let $X$ be a Banach space. We say that $x\in X$ is {\sl extreme} if it
is an extreme point of the unit ball of $X$ (in particular, $\|x\|=1$).
We say that $x\in\partial \text{Ball}(X)$ is {\sl exposed} in $X$ if
there exists a functional $L\in X^*$ that attains its norm at $x$ and
at no other point of the closed unit ball. The functional $L$ is often
assumed to be of norm 1 and is (then) called an {\sl exposing functional}
for $x$. Of course an exposed point is also extreme, but the converse
need not hold in general. However, we have the following lemma, the simple 
proof of which we omit. 
\begin{Lem} \label{lem:extexp}
   Let $X$ be a Banach space in which every point of unit norm is extreme.
   Then all points of unit norm are also exposed.
\end{Lem}

The concept of exposedness can be refined in the following manner.
We call $f\in \partial\text{Ball}(X)$ {\sl strongly exposed\/}
if there exists $L\in X^*$ with the properties: $L(f)=\|L\|=1$
and for any sequence $(f_n)_1^\infty$ in the ball of $X$
such that $\lim_{n\to\infty} L(f_n)=1$, it follows that $f_n$ converges
to $f$ in norm. It is not difficult to see that a strongly exposed point is
(indeed) exposed (by the same functional). By a theorem of Phelps \cite{phelps}, 
in a separable dual Banach space the closed unit ball is the closure of the 
convex hull of the strongly exposed points. In particular, because the Bergman space $A^1$ is
a dual space, section \ref{S:Bloch}, there exist many strongly 
exposed points.

\section{A criterion for strongly exposed points in $A^1$.} \label{S:SEPinA}
Let us illustrate these definitions for the Bergman space $A^1$.

\begin{Lem}
  Let $f\in A^1$ be of unit norm. Then the functional
$$
  L:g\in A^1 \mapsto \int_D g \overline{f}/|f| dA
$$
  is exposing for $f$. In particular, all functions of unit norm are exposed in 
  the unit ball of $A^1$.
  
\end{Lem}
It is not hard to show that all functions of unit norm in $A^1$ are extreme in the unit
ball. Thus the claim follows from Lemma \ref{lem:extexp}. However we opt for another 
proof:
\begin{proof}
  Let $f\in A^1$ be of unit norm. Suppose first that (some) $L\in (A^1)^*$
  is such that $L(f)=\|L\|=1$. Then, by the Hahn-Banach theorem, there
  exists $\psi\in L^\infty$ such that $\|\psi\|=1$ and
  $L(g)=\int_D g\overline\psi \,dA$ for all $g\in A^1$. Because also
  $\int_D |f| \,dA=1$ it follows that $\psi=f/|f|$
  (almost everywhere). In particular, an exposing functional for $f$
  is unique (if it exists).
  We finish the proof by showing that the functional
  $L(g)=\int_D g\overline{f} / |f|\, dA$ is indeed exposing for $f$.
  Clearly, $L$ attains its norm at $f$. Suppose that $g$ in the closed unit ball
  of $A^1$ is such that $L(g)=1$. Then by the above reasoning, $f/|f|=
  g/|g|$ almost everywhere on $D$. Hence $g/f$ is a positive
  meromorphic function on $D$, thus constant. Because $\|g\|=\|f\|=1$,
  $g$ must equal $f$. This concludes the proof.
\end{proof}

By contrast, not all functions of unit norm are strongly exposed in the unit
ball of $A^1$. This is contained in the following proposition. 

\begin{Prop}\label{prop:logniet}
   The function $$f(z)=\frac{cz^2}{(1-z)^2\log^2(1-z)},$$ 
where $c$ is normalizing so that $\|f\|=1$,
is not strongly exposed in the unit ball of $A^1$. 
\end{Prop}
\begin{proof}
   For $-2<\beta<0$, let $f_\beta(z) = c_\beta (1-z)^\beta$, where the constant $c_\beta>0$ is normalizing,
   i.e. $\|f_\beta\|=1$. Let $\varphi_\beta = f_\beta/|f_\beta|$ and $\varphi_{-2}=\frac{1-\overline{z} }{ 1-z}$. 
   By construction, $\int_D f_\beta \overline{\varphi_\beta}\, dA = 1$ 
   for all $\beta$. Let 
$$\varphi=f/|f|=\varphi_{-2}\cdot\frac{z\log(1-\overline{z})}{\overline{z}\log(1-z)}.$$ 
   Then the $f$-exposing functional $L$ is given by
$$
   L:g\in A^1 \mapsto \int_D g\overline{\varphi}\, dA.
$$
   As $\|\varphi_\beta-\varphi_{-2}\|_\infty \to 0$ for $\beta\downarrow -2$, it follows that 
   $\lim_{\beta\downarrow -2} \int_D f_\beta\overline{\varphi_{-2}} \,dA=1$. Next, because $\frac{z\log(1-
   \overline{z})}{\overline{z}\log(1-z)}\to 1$ as $D\ni z\to 1$, the bounded function $\varphi-\varphi_{-2}$ is continuous on 
   $\overline{D}\setminus\{0\}$ and vanishes at $z=1$. Because $\|f_\beta\|\to\infty$ as $\beta\downarrow -2$, the normalizing 
   constants $c_\beta$ tend to $0$ as $\beta\downarrow -2$, thus the functions $f_\beta$ tend to $0$ uniformly on 
   $\overline{D}\setminus B(1,\varepsilon)$ for every $\varepsilon>0$, as $\beta\downarrow -2$. Hence, $\int_D 
   f_\beta(\overline{\varphi}- \overline{\varphi_{-2}}) \,dA\to 0$ as $\beta\downarrow -2$. Consequently, 
$$L(f_\beta)=\int_D f_\beta\overline{\varphi}\, dA \to L(f)=1.$$
   Because the functions $f_\beta$ tend to zero pointwise, they do not converge to $f$ in norm. 
   This demonstrates that $f$ is not strongly exposed. 
\end{proof}

Now let
$$
 (A^1)^{\perp} = \{ \psi\in L^\infty : \int_D f\overline{\psi} \,dA = 0
  \text{ for all\,} f\in A^1\}
$$
denote the annihilator of $A^1$ contained in $L^\infty$.
The space $(A^1)^{\perp}$ is quite large but we do not know of a
structural description of its elements.
Finally let $C$ denote the continuous functions on $\overline{D}$.

We are now ready to give an abstract characterization of the strongly
exposed points of $A^1$.

\begin{Th} \label{th:jan}
  Let $f\in A^1$ be of unit norm. Then $f$ is strongly exposed if and
  only if the $L^\infty$-distance of $f/|f|$ to the space
  $(A^1)^{\perp} + C$ is less than one.
\end{Th}
\begin{proof}
  We first show that the distance condition is necessary.
  We argue by contradiction. Thus let $f\in A^1$ be strongly
  exposed and suppose the $L^\infty$-distance of $f/|f|$ to
  $(A^1)^{\perp} + C$ is one. Pick a point $z_0$ such that
  $f(z_0)\neq 0$. Let $A_{z_0}^1$ denote the subspace of
  all Bergman functions vanishing at $z_0$. We let $L'$ denote
  the restriction of the $f$-exposing functional $L:g\mapsto
  \int g\overline{f}/|f| \,dA$ to $A_{z_0}^1$. By the Hahn-Banach
  theorem the operator norm of $L'$ equals the $L^\infty$-distance
  of $f/|f|$ to $(A_{z_0}^1)^{\perp}$. Now if $\psi\in(A_{z_0}^1)^{\perp}$,
  then with the choice $c=\int_D \psi \,dA$, the function
  $\psi(w)-\frac{c}{ (1-\overline{z_0}w)^2}$ annihilates both
  $A_{z_0}^1$-functions and constants. This shows that
  $(A_{z_0}^1)^{\perp}\subset (A^1)^{\perp} + C$. By the
  assumption on $f/|f|$, we conclude that $L'$ has operator norm 1.
  Hence we find a sequence of functions $f_n$ in the unit ball of
  $A_{z_0}^1$ for which $L(f_n) = L'(f_n) \to 1$. Yet contrary to
  the assumption of strong exposedness of $f$, the functions $f_n$
  do not converge to $f$ in norm. Indeed, norm convergence implies
  pointwise convergence, which fails at the point $z_0$.

  Next we show that the distance condition is sufficient.
  This distance condition strongly resembles one in a theorem of
  the second author on the strongly exposed points in the Hardy
  space $H^1$ of the unit ball $B_n$ of ${\Bbb C}^n$~\cite{janw93}.
  There it is proven that an exposed point $F$ is strongly exposed in
  $H^1$ if and only if the $L^\infty$-distance of the function
  $F/|F|$ on the sphere $S$ of ${\Bbb C}^2$ to the space
  $(H^1)^{\perp} + C(S)$ is less than one.
  By Theorem 7.2.4 in~\cite{rudincn}, the Bergman space $A^1$ is
  isometrically contained in the Hardy space $H^1$ of the unit ball $B_2$
  in ${\Bbb C}^2$; namely, look at all holomorphic functions $F(z,w)$ on
  $B_2$ which depend only on $z$: $F(z,w)=F(z,0)$. Then $F$ is in
  $H^1(B_2)$ if and only if $f(z):=F(z,0)$ is in $A^1$ and the corresponding
  norms are then the same. Similarly, $(A^1)^\perp$ can be interpreted as
  a subspace of $(H^1)^\perp$ and $C(\overline{D})$ as a subspace of
  $C(S)$.

  Let our function $f\in A^1$ correspond with $F\in H^1$. Because
  $F$ is continuous on an open subset of $S$ (in fact, all of $S$
  except possibly ${\Bbb T} \times \{0\}$), $F$ is exposed in $H^1$.
  Because of the inclusion $(A^1)^\perp + C
  \subset (H^1)^\perp + C(S)$, $F/|F|$ has $L^\infty$-distance
  less than one to $(H^1)^\perp + C(S)$. Hence $F$ is strongly exposed in
  $H^1$ by~\cite{janw93}. Thus $F$, or rather $f$, is strongly exposed in
  $A^1\subset H^1$. This finishes the proof.
\end{proof}

The question now is: how can we estimate the distance in $L^\infty$
of $\varphi=f/|f|$ to $(A^1)^\perp + C$, where $f$ is a given function in
$A^1$? Clearly the distance cannot exceed one. Throughout the remainder
we will use various techniques to estimate said distances.

\section{The functions $(z-\alpha)^\beta$, Part I.} \label{S:polyn}
In order to simplify the necessary calculations we will test strong
exposedness on functions of a particularly simple form, i.e. polynomials.
It will later be shown (section~\ref{S:results}) that we may then even restrict to simple polynomials
of the form $f(z)=c(z-\alpha)^n$, where $c$ is normalizing. Having then
obtained our results for these functions it is easy to generalize to
functions of the form $f(z)=c(z-\alpha)^\beta$ for non-integer $\beta$
(in which case ${|\alpha|\ge 1}$, obviously).

So let us first look at polynomials: $f(z)=c(z-\alpha)^n$. We assume $n\ge 1$ because
unimodular constants are clearly strongly exposed. 
The case where $|\alpha|>1$ is the easiest: $f/|f|$ is
continuous on $\overline{D}$, so $f$ will be strongly exposed. When $|\alpha| < 1$
the proof that $f$ is again strongly exposed is a little more involved. Let us
write $\varphi=f/|f|$. If we can show that the Bergman projection $P\varphi$
is continuous on $\overline{D}$ (thus bounded), we will be done because
it will then follow that $\varphi=(\varphi-P\varphi) + P\varphi$ is
contained in $(A^1)^\perp + C$ rather than $(A^2)^\perp + C$.
Write $\varphi = \psi_1 + \psi_2$, where $\psi_1$ is compactly
supported in $D$ and $\varphi\equiv \psi_1$ on a neighborhood of
$\alpha$. From (\ref{eq:Bproj}) we see that $P\psi_1$ is holomorphic
across the unit circle because of the support of $\psi_1$. 
For the other function, using the series expansion (\ref{eq:bergP2}) for the Bergman
projection, we see that the smoothness of $\psi_2$ implies continuity (smoothness) 
of $P\psi_2$ on $\overline{D}$. This proves that $P\varphi$
is continuous on $\overline{D}$ and we conclude that $f$ is strongly
exposed.

In this section we can solve the case when $|\alpha|=1$ only partially, that is, depending
on the degree $n$ of the polynomial. We may then of course assume that
$\alpha=1$. Let us write $f_n(z)=c_n (1-z)^n$ and $\varphi_n = f_n/|f_n|$.
The corresponding exposing functional $L$ for $f_2$ is given by
$$
  L(g) = \int_D g(z) \frac{1-\overline{z} }{ 1-z} \,dA(z)
       = \int_D \frac{g(z) }{ 1-z} (1-\overline{z}) \,dA(z).
$$
Integrating first over circles we see that there exist constants
$C_0$ and $C_1$ (independent of $g$) such that
$L(g) = C_0 g(0) + C_1 g'(0)$. Thus there exists a polynomial
$\psi$ such that $L(g)=\int_D g\overline{\psi} \,dA$. (Alternatively,
verify that $\psi=P\varphi_2$ is a polynomial.) But this means
that $\varphi_2 - \psi$ is contained in the annihilator of $A^1$,
hence that $\varphi_2\in (A^1)^\perp + C$ and subsequently $f_2$
is strongly exposed.

Quite similarly one shows that for all {\sl even} $n$, $\varphi_n$ is
contained in $(A^1)^\perp + C$ and that $f_n$ is strongly exposed in $A^1$.

We come to the following ``odd" proposition on real powers.

\begin{Prop}\label{prop:sinest} Let $f_\beta(z)=c_\beta (1-z)^\beta$. Then for all $\beta>-1$, the $L^\infty$-distance of $\varphi_\beta$ to $(A^1)^\perp + C$ 
is at most $|\sin(\frac{\beta\pi}{2})|$. In particular, for all $\beta>-1$, $\beta\neq 1,3,5,\ldots$, the function 
$f_\beta$ is strongly exposed in the unit ball of $A^1$.
\end{Prop}
\begin{proof}
   Of course, there is nothing to prove for odd $\beta$, so we take $\beta>-1$ not odd. We will exploit the fact that the 
   functions $\varphi_0,\varphi_2,\varphi_4,\ldots$ are contained in the space $(A^1)^\perp + C$. We find an integer $n\ge0$ 
   such that $\beta\in (2n-1,2n+1)$. Let $\theta=|\beta-2n| < 1$. Because $\varphi_{a+b}=\varphi_a\varphi_b$, 
\begin{eqnarray}\nonumber
   \|\varphi_\beta - \cos(\frac{\pi\theta}{2})\varphi_{2n} \|_\infty &=& \| \varphi_\theta - \cos(\frac{\pi\theta}{2}) %
   \|_\infty =\sup_{|t|<\frac{\pi\theta}{2}}|e^{it} - \cos(\frac{\pi\theta}{2})  | \\   
   {}&=& \sin(\frac{\pi\theta}{2}) = |\sin(\frac{\pi\beta}{2})|.\nonumber
\end{eqnarray}
This gives the desired upper bound for the $L^\infty$-distance of $\varphi_\beta$ to $(A^1)^\perp + C$. 
By Theorem \ref{th:jan}, $f_\beta$ is strongly exposed. 
\end{proof}

In section~\ref{S:results} we investigate the odd powers in greater detail. 
Before doing so, we need to investigate the Bergman projection further.

\section{The Bloch space.} \label{S:Bloch}
Recall the Bergman projection $P:L^2\to A^2$,
$$
  Pf(z) = \int_D \frac{f(w)}{ (1-z\overline{w})^2} \,dA(w).
$$
We have already used the Bergman projection $P$ to prove strong
exposedness, namely in those cases where $P$ projects the {\sl
bounded} function $\varphi=f/|f|$ to a continuous function on $\overline D$.
However, a priori we cannot even expect $P$ to project bounded functions
to bounded functions. Obviously we would like to understand better how
$P$ acts on bounded functions. For this we need to discuss the Bloch space.

The {\sl Bloch space} $\cal B$ consists of all holomorphic functions $f$ on $D$
with the property that $(1-|z|^2) |f'(z)|$ is bounded on $D$.
Equipped with the norm
\begin{equation}
  \|f\|_{\cal B} := |f(0| + \sup_{z\in D} \, (1-|z|^2) |f'(z)|,\label{eq:blochnorm}
\end{equation}
$\cal B$ becomes a Banach space. The set of all functions $f$ in
$\cal B$ for which the expression $(1-|z|^2) |f'(z)| \to 0$ as $|z|\to 1$ is a closed
subspace of $\cal B$, called the {\sl little Bloch space} ${\cal B}_0$.
Finally, let $C_0 (D)$ denote the continuous functions on $\overline{D}$
that are zero on $\Bbb T$.

\begin{Th}[\cite{coiff}] \label{th:littleB}
  The Bergman projection $P$ maps $L^\infty$ boundedly onto $\cal B$.
  Furthermore, $P$ maps both $C$ and $C_0$ boundedly onto ${\cal B}_0$.
\end{Th}
\begin{proof}
  Cf.~\cite{hkz}, Theorem 1.12. 
\end{proof}

For future reference we remark that the proof of Theorem \ref{th:littleB} in \cite{hkz} gives that 
the norm of $P$ is at most $\frac{\pi}{8}$ and that if $f\in\mathcal{B}$ satisfies $f(0)=f'(0)=0$, then 
the $L^\infty$-function $\psi=(1-|w|^2) f'(w)/\overline{w}$ is mapped to $f$ under $P$.

It can be shown the Bloch space is the dual of the Bergman space $A^1$
,
 while the Bergman space is the dual of the little Bloch space $ \mathcal{B}_0$.
However, the resulting operator norms are equivalent with, but not equal to the standard norms that
 we defined previously. See \cite{hkz}, Chapter 1.
 The strongly exposed points
in the Bergman space under the operator norm have been described by C.~Nara~\cite{nara}. 

\section{The space $(A^1)^\perp + C$.} \label{S:annih}

We recall that $(A^1)^\perp +C$ plays the same role in Theorem~\ref{th:jan}
with respect to the Bergman space as $(H^1)^\perp +C({\Bbb T})$ does with respect
to $H^1$ of the unit ball in $\Bbb C$. In $\Bbb C$, the space $(H^1)^\perp +C$
is nothing other than the space $H^\infty + C({\Bbb T})$ which has been studied extensively. 
It is a famous result (\cite{helsar},\cite{rudinhc}) that 
$H^\infty\,+C({\Bbb T})$ is a {\sl closed} subspace of $L^\infty$. From this 
then it follows relatively easily that $H^\infty +C({\Bbb T})$ is in fact an 
algebra. We will now discuss how these results extend to the space 
$(A^1)^\perp + C$.

\begin{Th}\label{th:closed}
  The space $(A^1)^\perp +C$ is a proper, closed subspace of $L^\infty$.
\end{Th}
\begin{proof}
  The kernel of the map $P:L^\infty\to {\cal B}$ is $(A^1)^\perp$.
  Because ${\cal B}_0$ is closed in $\cal B$, $P^{-1}({\cal B}_0)$ is
  closed in $L^\infty$ by the continuity of $P$. By Theorem~\ref{th:littleB},
  $L^\infty\neq P^{-1}({\cal B}_0)=(A^1)^\perp +C$ and we are done.
\end{proof}

\begin{Th}\label{th:module}
  The space $(A^1)^\perp + C$ is a $C$-module.
\end{Th}

Before we give the proof we need a lemma.

Let $L_0^\infty$ be the subspace of $L^\infty$ consisting of all
$L^\infty$-functions that satisfy
$$
  \lim_{r\to 1} {\text{ess}\sup}_{r<|z|<1} |f(z)| = 0.
$$

\begin{Lem}\label{lem:zerob}
  The space $L_0^\infty$ is a closed subspace of $(A^1)^\perp + C$ and
  $(A^1)^\perp + C$ is closed under multiplication by functions in $L_0^\infty$.
\end{Lem}
\begin{proof}
  Clearly, $L_0^\infty$ is closed. Also, the product of a
  function in $L_0^\infty$ and a bounded function will again be in $L_0^\infty$
  so what remains is to show that $L_0^\infty$ is contained in $(A^1)^\perp + C$.
  Take $\psi\in L_0^\infty$. We write $\psi=\psi_1+\psi_2,$ where $\psi_1$ is the
  restriction of $\psi$ to the disc around zero with radius $r$. If $r$ is close
  enough to 1, then $\|\psi_2\|_\infty$ will be arbitrarily small by the
  assumption on $\psi$. Hence the ${\cal B}$-norm of $P\psi_2$ will be arbitrarily
  small by the continuity of $P$. On the other hand, $P\psi_1$ is holomorphic
  across ${\Bbb T}$, so $P\psi_1\in {\cal B}_0$. It follows that the ${\cal B}$-distance
  of $P\psi$ to ${\cal B}_0$ will be at most the ${\cal B}$-norm of $P\psi_2$,
  i.e. arbitrarily small. Because ${\cal B}_0$ is closed in ${\cal B}$,
  we conclude that $P\psi\in{\cal B}_0$. By the proof of Theorem~\ref{th:closed},
  $\psi\in (A^1)^\perp + C$ and we are done.
\end{proof}

We proceed with the proof of Theorem~\ref{th:module}:
\begin{proof}
  The space $(A^1)^\perp$ is closed under multiplication by $\overline{z}$.
  Because $(A^1)^\perp+C$ is closed and by the Stone-Weierstrass theorem,
  we need then only show that $zg(z)\in (A^1)^\perp + C$ when $g$ is in
  $(A^1)^\perp$.
  Take $f\in A_0^1=zA^1$, say $f(z)=zF(z), F\in A^1$. Then, with the
  $L^2$-inner product $\langle.,.\rangle$, $\langle f,zg\rangle=
  \langle F,|z|^2 g\rangle$. Observe that the function
%Because  $(1-|z|^2)g(z) \in L_0^\infty$,
  $|z|^2 g$ is contained in $(A^1)^\perp + C$ because $(1-|z|^2)g(z) \in L_0^\infty$ and by 
  Lemma~\ref{lem:zerob}. Let's say, $|z|^2 g = g_1 + \varphi_1$, where 
  $g_1\in (A^1)^\perp$ and $\varphi_1\in C$. Then
  $\langle f,zg\rangle=\langle F,\varphi_1\rangle $.
  Next we approximate $\varphi_1$ uniformly with a trigonometric polynomial
  $p_1=p_1(\varphi_1)$, i.e.~$\|\varphi_1-p_1\|_\infty < \varepsilon$.
  The integral $\langle F,p_1\rangle $ depends on the
  Taylor coefficients of $F$ in a finite fixed set of places. Because $f=zF$ has
  the same coefficients, albeit shifted, we can find a polynomial
  $p_2$ such that $\langle F,p_1\rangle =\langle f,p_2\rangle $, for all
  $f=zF\in A_0^1$.

  Now, $\langle f,zg\rangle =\langle F,\varphi_1\rangle =\langle F,p_1\rangle
  + \langle F,\varphi_1-p_1\rangle = \langle f,p_2 \rangle + \langle F,
  \varphi_1-p_1\rangle $, so that $|\langle f,zg-p_2\rangle| \le
  \varepsilon \|F\|_{A^1}$. We remark that the $A^1$-norms of $f$ and
  $F$ are equivalent in the sense that 
  for all $F\in A^1$: $$\|zF\|_{A^1} \le \|F\|_{A^1} \le 4
  \|zF\|_{A^1}.$$ Hence, $|<f,zg-P_2>| \le 4\varepsilon \|f\|_{A^1}$. 
  By the Hahn-Banach theorem, the $L^\infty$-distance of
  $zg-p_2$ to the annihilator of $A_0^1$ is at most $4\varepsilon$. And
  since $p_2$ is continuous, and $(A_0^1)^\perp + C = (A^1)^\perp + C$,
  the $L^\infty$-distance of $zg$ to $(A^1)^\perp + C$ is at most
  $4\varepsilon$, thus zero. By Theorem~\ref{th:closed},
  $zg\in (A^1)^\perp + C$ and the proof is complete.
\end{proof}

It is well-known that the space $H^\infty +C(\overline{D})$ is closed
in $L^\infty(D)$ (Theorem 6.5.5. in~\cite{rudincn}, \cite{rudinhc}). 
Let us write ${\cal A}:=L_0^\infty + \overline{H^\infty} + C({\overline D})$, where the bar denotes complex 
conjugation. By the preceding remarks, the space ${\cal A}$ is a 
non-trivial closed {\sl algebra} contained in $(A^1)^\perp + C$, and the
space $(A^1)^\perp + C$ is an ${\cal A}$-module. It should be stressed however
that $(A^1)^\perp + C$ is not an algebra.

\begin{Lem}
  Let $f_\beta = (1-z)^\beta$ and let $\varphi_\beta=f_\beta/|f_\beta|$
  for $\beta\in{\Bbb R}$. Then $\varphi_{-4}\in (A^1)^\perp$, but
  $\varphi_{-2}$ is not contained in the space $(A^1)^\perp +C$.
\end{Lem}
\begin{proof}
  Using the Stokes theorem one obtains that, at least formally, for every
  polynomial $F$:
\begin{eqnarray}\nonumber
  \int_D F \overline{\varphi_{-4}} \, dA &=& \int_D \overline{\partial} \Bigl[F(z) \frac{(1-z)^2}{(1-\overline{z})}\Bigr] \, dA \\
\nonumber  &=& \int_S F(z)\frac{ (1-z)^2 }{(1-\overline{z})}\, \frac{dz}{ 2\pi}
\\
\nonumber  &=& \int_S -F(z)(1-z)z \,\frac{dz}{ 2\pi} = 0.
\end{eqnarray}
(In fact, by the same argument, $\int_D F(z)\overline{\varphi_{-2k}} \,dA=0$
for all $k=2,3,4,\ldots$.) Here we have used the identity $z\overline{z}=1$
on $S$ to simplify the integrals over the circle. We conclude that (formally) 
$P\varphi_{-4}=0$, that is, $\varphi_{-4}\in (A^1)^\perp$. 
This claim can be made precise by a limit argument involving integration over the
unit disc with a small disc around the point $z=1$ punched out. Alternatively, one can 
directly calculate the Bergman projection of $\varphi_{-4}$. We omit the details. By Theorem \ref{th:jan}
combined with the calculations in the proof of Proposition \ref{prop:logniet}, we conclude that the $L^\infty$-distance
of $\varphi_{-2}$ to $(A^1)^\perp + C$ is $1$. In particular, the function $\varphi_{-2}$ is not contained in
$(A^1)^\perp + C$. 
\end{proof}

\begin{Cor}
The space $(A^1)^\perp + C$ is not an algebra.
\end{Cor}
\noindent Indeed, $\varphi_{-4}$ and $\varphi_2$ are both contained in $(A^1)^\perp + C$,
but their product $\varphi_{-4}\cdot \varphi_2 = \varphi_{-2}$ is not.

Next, let $u$ be an automorphism (M\"{o}bius map) of $D$. If $\psi$ is an
element of $L^\infty(D)$ one can define the composition $\psi\circ u$ in
$L^\infty$ of $\psi$ and $u$ as (represented by) the composition
of $\Psi$ with $u$, where $\Psi$ is any representative of $\psi$. That this
yields a well-defined element of $L^\infty$ follows from the fact that
$u$ and its inverse map sets of Lebesgue measure zero to sets of
Lebesgue measure zero. It is easily seen that the map $\psi\mapsto
\psi\circ u$ is an isometric isomorphism of $L^\infty$.

\begin{Prop}\label{prop:autom}
  The space $(A^1)^\perp + C$ is invariant under composition with
  automorphism of $D$.
\end{Prop}
\begin{proof}
  Clearly, the space $C$ is invariant under composition with automorphisms of $D$.
  Take an element $g\in (A^1)^\perp$, and let $u$ be an automorphism of
  $D$. We will show that $g\circ u$ is contained in $(A^1)^\perp +C$.
  Let $f$ be an element of $A^1$. Then $\int_D f \overline{g\circ u} \, dA=
  \int_D (f\circ u^{-1}) \overline{g}  J_{\Bbb R} (u^{-1}) \, dA$, where
  $J_{\Bbb R} (u^{-1})$ is the real Jacobian of $u^{-1}$, an element of $C$.
  By Theorem~\ref{th:module} there exist $g^*\in (A^1)^\perp$
  and  $h\in C$ such that $g J_{\Bbb R} (u^{-1}) = g^* + h$.
  Thus, because $f\circ u^{-1}$ is contained in $A^1$, $\int_D f \overline{g\circ u}
  \, dA = \int_D (f\circ u^{-1}) \overline{h} \, dA = \int_D
  f \overline{(h\circ u)} J_{\Bbb R}(u) \, dA$. We conclude that $g\circ u -
  (h\circ u) J_{\Bbb R}(u)$ annihilates the Bergman space, hence
  $g\circ u\in (A^1)^\perp + C$.
\end{proof}

\begin{Prop}\label{prop:newsep}
  Let $f$ be a strongly exposed point in $A^1$.
\begin{itemize}
\item[(a)] If $u$ is an automorphism of $D$, then the normalized function
           $f_1=C_1(f\circ u)$ is strongly exposed.
\item[(b)] If $v\in C$ is holomorphic on $D$ and zero-free on the circle,
           then the normalized function $f_2=C_2 f v$ is strongly
           exposed.
\end{itemize}
  Furthermore, the functions $\varphi=f/|f|, \varphi_1=f_1/|f_1|$ and
  $\varphi_2=f_2/|f_2|$ have the same $L^\infty$-distance to $(A^1)^\perp + C$.
\end{Prop}
\begin{proof}
  $(a)$\quad There exist $g\in (A^1)^\perp, h\in C$ such that $\|\varphi-
  g-h\|_\infty < 1$. By Proposition~\ref{prop:autom} $g\circ u$ is again
  contained in $(A^1)^\perp +C$. Because $\varphi_1=\varphi\circ u$:
  $\| \varphi_1 - g\circ u -h\circ u\|_\infty<1$, and we conclude that $f_1$
  is strongly exposed. Also, the $L^\infty$-distance of $\varphi_1$
  to $(A^1)^\perp + C$ does not exceed that of $\varphi$. Replacing
  $u$ by its inverse, the reverse inequality follows. %\linebreak
  $(b)$\quad With $g$ and $h$ as above and $\varphi_2=\varphi \frac{v}{|v|}$,
  $\|\varphi_2 - g\frac{v}{ |v|} - h\frac{v}{|v|}\|_\infty < 1$.
One finishes the proof as before, using Lemma~\ref{lem:zerob} and the
fact that $\frac{v}{|v|}$
  is invertible in $L_0^\infty + C$.
\end{proof}

\section{The functions $(z-\alpha)^\beta$, Part II.} \label{S:results}  
We saw in section \ref{S:SEPinA} that the functions $f_\beta=c_\beta (1-z)^\beta$ are strongly
exposed in the unit ball of $A^1$ for all $\beta>-1$ except {\sl possibly} when $\beta=1,3,5,\ldots.$ This was
deduced from rather straightforward estimates of the $L^\infty$-distances of the functions $\varphi_\beta=f_\beta/|f_\beta|$
to the space $(A^1)^\perp+C$ (Proposition \ref{prop:sinest}). In this section we will sharpen these estimates
and answer the question of strong exposedness for odd exponents.  

\begin{Th}{\label{th:betas}}   
  For all $\beta\ge 0$, the Bloch distance of the function $P\varphi_\beta$ 
  to ${\cal B}_0$ equals $\frac{4}{\pi} \frac{|\sin(\frac{\beta\pi}{2})|}{\beta+2}$.
\end{Th} 
\begin{proof} 
  We showed in section~\ref{S:polyn} that the functions $P\varphi_{2n}$  are contained in $\mathcal{B}_0$ so 
  henceforth we will assume that $\beta$ is not even. It is convenient to rewrite $\varphi_\beta$ as  $\varphi_\beta (w) 
  = (1-w)^{\beta/2} / (1-\overline{w})^{\beta/2}$. Using the series expansions for the Bergman kernel $1/(1-z\overline{w})^2$
  (see (\ref{eq:bergP2})), as well as for $(1-w)^{\beta/2}$, and $1/(1-\overline{w})^{\beta/2}$, we evaluate the Bergman 
  projection $P\varphi_\beta$. One obtains $P\varphi_\beta=\sum_{n=0}^\infty c_{\beta,n} z^n$, where 
$$c_{\beta,n} = \frac{n+1}{ 
  \Gamma(-\frac{\beta}{2}) \Gamma(\frac{\beta}{2})} \sum_{m=0}^\infty \frac{\Gamma(m+\frac{\beta}{2}) 
  \Gamma(m+n-\frac{\beta}{2})}{m! (m+n+1)!}. 
$$   
  we claim that for fixed $\beta>0$:
\begin{equation} 
  \sum_{m=0}^\infty \frac{\Gamma(m+\frac{\beta}{2}) \Gamma(m+n-\frac{\beta}{2})}{m! (m+n+1)!} 
   = \frac{4}{n^2\beta(\beta+2)} (1+o(1)), \label{eq:size} 
\end{equation}   
  where the $o(1)$-term tends to zero as $n\to\infty$. This implies that  
$$ 
  c_{\beta,n} = \frac{4}{\Gamma(-\frac{\beta}{2})\Gamma(\frac{\beta}{2})}              
   \frac{1}{n\beta(\beta+2)}  (1+o(1))=
  \frac{-2\sin(\frac{\beta\pi}{2})}{\pi(\beta +2)n}(1+o(1)),
$$ 
  where the $o(1)$-term vanishes as $n\to\infty$. (Here we have used the functional equations $\Gamma(z+1)=z\Gamma(z)$
  and $\Gamma(z)\Gamma(1-z)\sin(\pi z)=\pi$.) But then,   
$$\lim_{x\uparrow 1} |(1-x^2) (P\varphi)'(x)| = 
  \frac{4|\sin(\frac{\beta\pi}{2})|}{\pi(\beta+2)},
$$
  so the Bloch distance of $P\varphi_\beta$ to ${\mathcal B}_0$ is at least $ \frac{4|\sin(\frac{\beta\pi}{2})
  |}{\pi(\beta+2)}$. On the other hand, for large $N$,   
$$
|(\sum_{n=N}^\infty c_{\beta,n}z^n)'| \le 
  \sum_{n=N}^\infty n|c_{\beta,n}| |z|^{n-1} 
  \le \frac{2|\sin(\frac{\beta\pi}{2})|}{\pi (\beta+2)}\cdot \frac{1+o(1)}{1-|z|}, 
$$ 
  where the $o(1)$-term tends to zero as $N$ increases. Using the fact that the polynomials are contained in 
  ${\mathcal B}_0$ it follows that the Bloch distance of $P\varphi_\beta$ to ${\mathcal B}_0$ is at most 
  $\frac{4|\sin(\frac{\beta\pi}{2})|}{\pi(\beta+2)}$.
This then proves the theorem.  

\bigskip
  We turn to the claim (\ref{eq:size}). Let us first assume $\beta>2$. Given any large $n\in{\Bbb N}$, let $M=M_n$ be
  the integer nearest to $\sqrt{n}$. We write  
$$\sum_{m=0}^\infty \frac{\Gamma(m+\frac{\beta}{2}) \Gamma(m+n-\frac{\beta}{2})}{ 
  m! (m+n+1)!}=\sum_{m=0}^{M-1} + \sum_{m=M}^\infty.
$$ 
  Because $\beta > 2$, $\frac{\Gamma(m+\frac{\beta}{2})}{m!}$ is increasing in $m$. 
  On the other hand, $\frac{\Gamma(n+m-\frac{\beta}{2})}{(n+m+1)!}$ is decreasing in $m+n$. The first sum can thus 
  be estimated by 
$$\sum_{m=0}^{M-1} \le M \frac{\Gamma(M+\frac{\beta}{2})}{(M)!} 
  \frac{\Gamma(n-\frac{\beta}{2})}{(n+1)!}. 
$$ 
  Recall Stirling's formula:
\begin{eqnarray}\label{eq:stirling}
   \lim_{x\to\infty} \frac{\Gamma(x+1)}{\sqrt{2\pi x}} (\frac{x}{e})^x = 1.
\end{eqnarray}
  By this result, there exists a constant $A=A_\beta$, independent of $n$, such that  
$$\sum_{m=0}^{M-1} \le A \frac{M\cdot M^{\frac{\beta}{2}-1}}{ 
  n^{2+\frac{\beta}{2}} } = \frac{A}{n^2} \Bigl(\frac{M}{n}\Bigr)^{\frac{\beta}{2}}. 
$$
  Hence
\begin{equation}
  \sum_{m=0}^{M-1} = \frac{o(1)}{n^2}, \label{eq:first}
\end{equation}
  as $n\to\infty$. In the remaining sum, $\sum_{m=M}^\infty$, all the arguments in the Gamma functions and factorials 
  tend to infinity as $n\to \infty$. Another application of Stirling's formula seems in place. One obtains that, given any
  $\varepsilon>0$, for all sufficiently large $n$ and all $m\ge M$,

$$   
  \Big|\, \frac{\Gamma(m+\frac{\beta}{2}) \Gamma(m+n-\frac{\beta}{2})}{ 
  m! (m+n+1)!}  \quad \Big/ \quad \frac{m^{\frac{\beta}{2}-1}}{(n+m)^{\frac{\beta}{2}+2}} \quad - 1 \,\Big| < \varepsilon.
$$  
  In particular, 
$$
  \Big|\, \sum_{m=M}^\infty \frac{\Gamma(m+\frac{\beta}{2}) \Gamma(m+n-\frac{\beta}{2})}{ 
  m! (m+n+1)!} \quad \Big/ \quad \sum_{m=M}^\infty \frac{m^{\frac{\beta}{2}-1}}{ 
  (n+m)^{\frac{\beta}{2}+2}} \quad -1 \Big| <\varepsilon, 
$$  
  as $n\to\infty$. Therefore, by (\ref{eq:first}), the claim (\ref{eq:size}) follows once we show that  
$$
  \sum_{m=M}^\infty \frac{m^{\frac{\beta}{2}-1}}{ 
  (n+m)^{\frac{\beta}{2}+2}}\quad \Big/ \quad \frac{4}{n^2 \beta(\beta+2)} \quad \to\quad 1, 
$$
  as $n\to\infty$. 
  Let us investigate the functions $g_n(x) = \frac{x^{\frac{\beta}{2} -1}}{(n+x)^{\frac{\beta}{2}+2}}$.
  For all $x\ge 1$, $g_n(x)\le \frac{1}{x(n+x)^2}$, so $g_n(x)\le \frac{1}{n^{5/2}}$ when $x\ge M$. 
  There is a number $x_{\beta,n} > 0$ such that $g_n(x)$ is increasing on the interval $(0, x_{\beta,n}]$ and decreasing 
  on the interval $[x_{\beta,n},\infty)$. Hence, the sum $\sum_{m=M}^\infty \frac{m^{\frac{\beta}{2}-1}}{(n+m)^{\frac{\beta}{2}+2}}$ 
  and the integral $\int_{M}^\infty \frac{x^{\frac{\beta}{2}-1}}{(n+x)^{\frac{\beta}{2}+2} }
  \,dx$ differ at most $\frac{4}{n^{5/2}}=\frac{o(1)}{n^2}$. By a change of variables,  
$$
  \int_{M}^\infty  \frac{x^{\frac{\beta}{2}-1}}{(n+x)^{\frac{\beta}{2}+2} }\,dx=\frac{1}{n^2} 
  \int_{\frac{M}{n}}^\infty \frac{x^{ \frac{\beta}{2}-1 }}{(1+x)^{ \frac{\beta}{2}+2}}\, dx. 
$$ 
  Now, with $B(.,.)$ the standard Beta-function, 
$$
  \int_0^\infty \frac{x^{\frac{\beta}{2}-1}}{(1+x)^{\frac{\beta}{2}+2} }\,dx= B(\frac{\beta}{2},2) = \frac{4}{\beta(\beta+2)}. 
$$ 
  On the other hand, as $n\to\infty$,  
$$
  \int_0^{\frac{M}{n}} \frac{x^{\frac{\beta}{2}-1}}{(1+x)^{\frac{\beta}{2}+2}} \,dx = o(1). 
$$ 
  By the preceding estimates, the claim (\ref{eq:size}) now follows for all $\beta>2$. 

\quad  When $0<\beta<2$ we proceed as follows. Given a large $n\in{\Bbb N}$, we let $M=M_n$ be the integer nearest to $n^{\frac{\beta}{4}}$.
  Now the terms in the sum 
$$\sum_{m=0}^\infty \frac{\Gamma(m+\frac{\beta}{2}) \Gamma(m+n-\frac{\beta}{2} )}{ 
  m! (m+n+1)!}=\sum_{m=0}^{M-1} + \sum_{m=M}^\infty.
$$ 
  are decreasing. The first sum can be estimated by 
$$
  \sum_{m=0}^{M-1} \le M \Gamma(\frac{\beta}{2})\frac{\Gamma(n-\frac{\beta}{2})}{(n+1)!}\le 
  \frac{A_\beta}{n^{2+\frac{\beta}{4}}}=\frac{o(1)}{n^2}.
$$
  The second sum can be dealt with as before. (Now the functions $g_n(x)$ are decreasing on $(0,\infty)$, which makes the
  analysis even simpler.) We omit the details. This finishes the proof of equation (\ref{eq:size}) for all $\beta>0$.
\end{proof}   

\begin{Cor}
   Let $d(\varphi_\beta, (A^1)^\perp + C)$ denote the $L^\infty$-distance of $\varphi_\beta$ to
  $(A^1)^\perp +C$. Then for all $\beta\ge 0$,
\begin{equation}
  \frac{1}{2} \frac{|\sin(\frac{\beta\pi}{2})|}{\beta+2} \le
  d(\varphi_\beta, (A^1)^\perp + C) \le 
  \frac{4}{\pi} \frac{|\sin(\frac{\beta\pi}{2})|}{\beta+2}\le\frac{2}{\pi},\label{eq:AenB}
\end{equation}
  In particular, all $f_\beta$ are strongly exposed for $\beta\ge 0$. 
 \end{Cor}

\begin{proof}
  Let $q:{\mathcal B}\to{{\mathcal B}\big/{\mathcal B}_0}$ be the quotient map.
  By Theorem~\ref{th:littleB}, the map $q\circ P: L^\infty\to {{\mathcal B}\big/{\mathcal B}_0}$ is 
  continuous and surjective. In the proof of Theorem~\ref{th:closed} it was shown that the kernel of the map 
  $q\circ P$ is the space $(A^1)^\perp + C$. It follows that the derived map
$$
  P^*:L^\infty \Big/ \Big((A^1)^\perp+C\Big) \to  {\mathcal B}\Big/{\mathcal B}_0
$$
  is bijective and bounded by $\frac{8}{\pi}$ (cf.~the proof of Theorem \ref{th:littleB}).  
  This gives the lower bound for $d(\varphi_\beta, (A^1)^\perp + C)$, because $\|P^*\varphi_\beta\|
  =\frac{4}{\pi} |\sin(\frac{\beta\pi}{2})| /\beta+2$. 
  
  By the closed graph theorem, the inverse ${P^*}^{-1}$ of $P^*$ is also bounded. Actually, we will show directly that  $\|{P^*}^{-1}\| \le 1$, 
  which in turn yields the upper bound for $d(\varphi_\beta, (A^1)^\perp + C)$. Let us suppose that $F\in {\mathcal B}/{\mathcal B}_0$ has 
  norm $1$. We need to show that ${P^*}^{-1}(F)$ has norm at most $1$ in $L^\infty\big/ ((A^1)^\perp+C)$. For any $\varepsilon>0$,
  we can find a representative $f\in\mathcal{B}$ of the coset $F$ such that $\|f\|_{\mathcal{B}} < 1+\varepsilon$. 
  We recall from the proof of Theorem \ref{th:littleB} that 
$$
   \psi(w) = (1-|w|^2)\cdot \frac{f'(w)-f'(0)}{\overline{w}}\in L^\infty
$$
  satisfies $f(z)-P\psi(z) = f(0) + f'(0)z\in{\mathcal B}_0$. Thus $\psi$ is a representative of ${P^*}^{-1}(F)$ in $L^\infty$. 
  Hence, by Lemma \ref{lem:zerob}, 
\begin{eqnarray}\nonumber
   \|{P^*}^{-1}(F)\|_{L^\infty/((A^1)^\perp+C)} &\le& d(\psi, (A^1)^\perp + C)\le \lim_{r\to 1} {\text{ess}\sup}_{r<|w|<1} |\psi(w)| \\
    &=& \limsup_{|w|\to 1} |(1-|w|^2)f'(w)| \le \|f\|_{\mathcal B} < 1+\varepsilon.\nonumber
\end{eqnarray}
\end{proof}
\begin{Cor}\label{cor:factors}
   Suppose that $g\in H(D)\cap C$ vanishes nowhere on $\Bbb T$. Let $z_1, z_2,\ldots,z_n\in{\Bbb T}$ be distinct and let
   $\beta_1,\beta_2,\ldots,\beta_n$ be real numbers greater than $-2$. Then the normalized function 
   $f(z)=cg(z)\prod_{i=1}^n (1-z\overline{z_i})^{\beta_i}$ is strongly exposed in the unit ball of $A^1$ if and
   only if all functions $f_{\beta_i}=c_{\beta_i}(1-z)^{\beta_i}$ are strongly exposed. In particular, 
   all choices of $\beta_i>-1$ yield strongly exposed points and all normalized polynomials are strongly exposed 
   in the unit ball of $A^1$.
\end{Cor}
\begin{proof}
   By part (b) of Proposition \ref{prop:newsep}, the factor $g(z)$ has no effect on strong exposedness of
   the function $f$.  
   Let $d_i=d(\varphi_{\beta_i},(A^1)^\perp + C)$ and let $\varphi=f/|f|$. We will show that
   $d(\varphi,(A^1)^\perp + C)=\max_i d_i$, which will give the desired result. 

   We find small pairwise disjoint neighborhoods $U_i$ of the $z_i$ and a partition
   $\chi_i$ of the unity relative to the $U_i$'s and $\overline{D}$. That is to say, we find
   continuous functions $\chi_i\ge 0$ on $\overline{D}$ such that $\chi_i\equiv 1$ on
   $U_i$ and $\sum_i \chi_i \equiv 1$ on $\overline{D}$. Then $\varphi=\sum_i \chi_i\varphi
   =\sum_i \varphi^{(i)}$. For every $i$, there exists a unimodular constant $\lambda=\lambda_i$ 
   for which $\varphi^{(i)}(z) -\lambda\varphi_{\beta_i}(z\overline{z_i})\in C$. Consequently,
   $d_i=d(\varphi^{(i)},(A^1)^\perp + C)$. Using the $C$-module structure of $(A^1)^\perp + C$ and the
   fact that $\chi_i \le 1$, it is easily seen that $d(\varphi^{(i)},(A^1)^\perp + C)\le
   d(\varphi,(A^1)^\perp + C)$, hence $\max_i d_i \le d(\varphi,(A^1)^\perp + C)$. 
   Conversely, if the functions $g_i\in (A^1)^\perp + C$ are such that $\|\varphi^{(i)} - g_i\|_\infty 
   < d_i +\varepsilon$, then 
$$\| \sum_i \chi_i\varphi^{(i)} - 
   \sum_i \chi_i g_i\|_\infty < \max_i d_i +\varepsilon.
$$ 
   Because $\varphi-\sum_i \chi_i\varphi^{(i)}=\sum_i \chi_i(1-\chi_i)\varphi\in C$ and
   $\chi_ig_i\in (A^1)^\perp + C$, it follows that $d(\varphi,(A^1)^\perp + C) <
   \max_i d_i + \varepsilon$. 
\end{proof}

We will now show that an estimate analogous to inequality (\ref{eq:AenB}) also holds for $\beta\in (-2,0)$. 

\begin{Lem}\label{lem:g_n}
   Let $g_n = \frac{1}{n} (1-z)^{-2+\frac{1}{n}}$. Then $\lim_{n\to\infty} \|g_n\|_1 = 1$. 
\end{Lem}
\begin{proof}
   One can perform the calculation
$$
   \int_D |1-z|^{-2+\frac{1}{n}} \,dA(z) = \sum_{k=0}^\infty \frac{\Gamma^2(k+1-\frac{1}{2n})}{
   \Gamma^2(1-\frac{1}{2n}) k! (k+1)!}.
$$
Using polar coordinates for the integral or hypergeometric functions for the sum, this expression can
be evaluted to ${\Gamma(1/n)}\big/{\big(\Gamma(1+1/2n)\big)^2}$.
 However the following asymptotics are fairly simple and useful in the following
Proposition.
 
   The terms $g_{n,k} = \frac{\Gamma^2(k+1-\frac{1}{2n} )}{k! (k+1)!}$ are decreasing in $k$.
   For large $n$, let $K=K_n=[\sqrt{n}]$. Then $\sum_{k=1}^K g_{n,k} \le \sqrt{n}$. In the remaining sum, 
   we can approximate the terms using Stirling's formula (\ref{eq:stirling}): $g_{n,k} \sim \frac{1}{k^{1+\frac{1}{n}}}$.
   Therefore, $\sum_{k=K+1}^\infty g_{n,k} \sim \int_{\sqrt{n}}^\infty \frac{1}{x^{1+\frac{1}{n}} }\sim n$. This proves the claim. 
\end{proof}

\begin{Prop}
   For all $-2<\beta<0$,
\begin{equation}
    \frac{2}{\pi} \frac{|\sin(\frac{\beta\pi}{2})|}{\beta+2}\le d(\varphi_\beta, (A^1)^\perp+C),\label{eq:negbeta}
\end{equation}
   where again $d(\varphi_\beta, (A^1)^\perp + C)$ denotes the $L^\infty$-distance of $\varphi_\beta$ to
   $(A^1)^\perp +C$.  
\end{Prop}
\begin{proof}
   Fix any $\beta\in (-2,0)$ and let $L=L_\beta$ be the functional $L:g\in A^1\mapsto \int_D g\overline{\varphi_\beta} \,dA$. We will show 
   that for the sequence $g_n$ from Lemma \ref{lem:g_n}, 
\begin{equation}
   \lim_{n\to\infty} L(g_n) = \frac{2}{\pi} \frac{|\sin(\frac{\beta\pi}{2})|}{\beta+2}.\label{eq:beta-2}
\end{equation}
   From this we will get the desired lower bound as follows. The functions $g_n$ and all their
   derivatives tend to zero uniformly on compact subsets of $D$. Given any (fixed) integer $N$ 
   we define functions $g_n^*(z):= g_n(z) - \sum_{k=0}^{N-1} \frac{g_{n}^{(k)}(0)}{k!} z^k$. It follows that
   $\lim_{n\to\infty} \|g_n^*\|_1 =1$ and 
\begin{equation}
   \lim_{n\to\infty} L(g_n^*) =\frac{2}{\pi}\frac{|\sin(\frac{\beta\pi}{2})|}{\beta+2}.\label{eq:beta-2*}
\end{equation}
   Furthermore, by construction, the first $N$ derivatives of the $g_n^*$ vanish at the origin. 
   If we let $z^N\! A^1\subset A^1$ denote the closed subspace of all functions in $A^1$ whose first $N$ derivatives vanish at
   the origin, then by equation (\ref{eq:beta-2*}), the norm of the functional $L$ restricted to $z^N\! A^1$ is at least $\frac{2}{\pi}
   |\sin(\frac{\beta\pi}{2})|/(\beta+2)$. By the Hahn-Banach theorem, the $L^\infty$-distance of $\varphi_\beta$ to $(z^N\! A^1)^\perp$ 
   is at least $\frac{2}{\pi}|\sin(\frac{\beta\pi}{2})|/(\beta+2)$. Consequently, the $L^\infty$-distance of $\varphi$ to ${\mathcal P} = 
   \cup_{N=1}^\infty (z^N\! A^1)^\perp$ is at least $\frac{2}{\pi}|\sin(\frac{\beta\pi}{2})| / (\beta+2)$. Observe that ${\mathcal P}$
   is uniformly dense in $(A^1)^\perp + C$, because it contains $(A^1)^\perp$ and all trigonometric polynomials. Therefore,
   the $L^\infty$-distance of $\varphi_\beta$ to $(A^1)^\perp + C$ is at least $\frac{2}{\pi}|\sin(\frac{\beta\pi}{2})| / (\beta+2)$.
     
\quad Let us turn to formula (\ref{eq:beta-2}). We calculate 
$$
   n L(g_n)=\int_D \frac{1}{(1-z)^{2+\frac{\beta}{2}-\frac{1}{n}}}
   \frac{1}{(1-\overline{z})^{-\frac{\beta}{2}}} \, dA(z)
$$ 
   using the series expansions for $(1-z)^\alpha$ and
   $(1-\overline{z})^\alpha$. After a routine calculation, one obtains the following expression:
$$
   n L(g_n)=\frac{1}{\Gamma(2+\frac{\beta}{2}-\frac{1}{n})\Gamma(-\frac{\beta}{2})} 
   \sum_{k=0}^\infty \frac{\Gamma(2+\frac{\beta}{2}-\frac{1}{n}+k)}{(k+1)!} \frac{\Gamma(k-\frac{\beta}{2})}{k!}.
$$
   Now 
$$    
   \frac{1}{\Gamma(2+\frac{\beta}{2}-\frac{1}{n})\Gamma(-\frac{\beta}{2})}\quad\to\quad
   \frac{1}{\Gamma(2+\frac{\beta}{2})\Gamma(-\frac{\beta}{2})}=
   \frac{2}{\pi}\frac{|\sin(\frac{\beta\pi}{2})|}{\beta+2},
$$
   as $n\to\infty$. So what's left to do is to show that 
$$
   \frac{1}{n}\sum_{k=0}^\infty \frac{\Gamma(2+\frac{\beta}{2}-\frac{1}{n}+k)}{(k+1)!} \frac{\Gamma(k-\frac{\beta}{2})}{k!}\quad\to\quad 1,
$$   
   as $n\to\infty$. This can easily be done by following the proof of Lemma \ref{lem:g_n}. 
\end{proof}

We conclude with a conjecture on the functions $f_\beta$ for $-2<\beta<0$, for which
strong exposedness is already implied by Proposition \ref{prop:sinest} when $-1<\beta<0$. 
Note that inequality (\ref{eq:negbeta}) is ``asymptotically sharp'' for $\beta\downarrow -2$:
$$
   \lim_{\beta\downarrow -2} \frac{2}{\pi} \frac{|\sin(\frac{\beta\pi}{2})|}{\beta+2} = 1 =
   d(\varphi_{-2},(A^1)^\perp + C).
$$ 

\begin{Conj}
   For all $-2<\beta<0$, $d(\varphi_\beta,(A^1)^\perp+C)=\frac{2}{\pi} \frac{|\sin(\frac{\beta\pi}{2})|}{\beta+2}$.
   In particular, the functions $f_\beta$ are strongly exposed for all said $\beta$. 
\end{Conj}

\end{document}